\documentclass{amsart}
\usepackage{amsmath}
\usepackage[latin1]{inputenc}
\usepackage[english]{}
\usepackage{graphicx}
\font\tengoth=eufm10
\font\sevengoth=eufm7
\font\fivegoth=eufm5
\newfam\gothfam
\textfont\gothfam=\tengoth \scriptfont\gothfam=\sevengoth
\scriptscriptfont\gothfam=\fivegoth

\def\blacksquare{\hbox to .60em{\vrule width .60em height .60em}}

\catcode`\@=12  
\def\C {\mathcal }

\def\hb {\hfil \break}
\def\n {\vskip 0.2cm \noindent }
\def\scirc{\,{\raise 0.8pt\hbox{$\scriptstyle\circ$}}\,}
\def\ins{\,{\raise 0.2cm \hbox{ $\scriptstyle \circ$}}\,}

\def\s{\sigma}

\def\h{\eta}

\author{Daniel Lehmann}
\address{4 rue Becagrun, 30980 Saint Dionisy, France}
\email{lehm.dan@gmail.com}
\title[Ordinary differential operators and curvilinear webs]{Ordinary and calibrated differential operators\hb 
Application to curvilinear webs} 
\date{}

\begin{document}
	\abstract{}{We study the space of the   solutions $s$ of any  system of partial differential equations 
		$D(j^ks)=0$ defined by a linear and homogeneous differential operator $D:J^kE\to F$ of any order $k\geq 1$,  which is ``ordinary" (i.e. which is generic in some sense among all $D$'s), $E$ and $F$ being vector bundles over a $n$-dimensional manifold $V$, and $D$ being assumed to be surjective at any point of $V$. In some range of  the ranks $p$ and $q$  of these bundles ($p<q\leq np$ in the case  $k=1$), we first give an upper-bound $\pi(n,k,p,q)$ for the dimension of the space ${\mathcal S}_m$ of the germs of solutions at a generic point $m$ of the ambiant manifold. If these ranks satisfy moreover to some condition of integrality (in the case  $k=1$, $\frac{p(n-1)}{q-p}$ must be  an integer), and we then say that $D$ is ``calibrated",  we build  a vector bundle $\mathcal E$ of rank $\pi(n,k,p,q)$ on $V$,  provided with a tautological connection $\nabla$,  whose curvature is  an obstruction for the dimension of ${\mathcal S}_m$ to reach its maximal value. We also prove a ``theorem of concentration'' : relatively to some convenient trivialization of $\mathcal E$,  some coefficients of this  curvature vanish systematically. 
		
		As an example, we provide, for any curvilinear  $d$-web on $V$, a differential operator $D$ of order one, which is always ordinary and calibrated, and for which ${\mathcal S}_m$ is the space of germs of abelian relations ([L]). Thus, we recover the Damiano's upper-bound ([D1]) for the rank of such a web, and we can define in the most general case    the  ``curvature'' of such a web, already  known for $n=2$  (see [BB] if $d=3$, and [Pa],[H1],[Pi1] for arbitrary  $d$), obstruction for this rank to be maximum.

		  }}
		\maketitle
			
	\section{Preliminaries}
		The framework is real $C^\infty$ or holomorphic. Exceptionally,     for theorem 1 in the    ``non-calibrated'' case,   analyticity (real or complex)  will be required.

		We are given two vector bundles  $E$ and  $F$    over  a  $n$-dimensional manifold $V$, and  a linear  homogeneous differential operator of   order $k$, ($k\geq 1$),  i.e.  a linear morphism of vector bundles   $D: J^kE\to F $, where $J^kE$ denotes the vector bundle of the $k$-jets of section of $E$.

		More generally,  for any integer $h \ (\geq k)$, we define by successive derivations  the $(h-k)^{th}$ prolongation   $D_h:J^{h}E\to J^{h-k}F $ of  $D$.  
		
		 Denoting  by $S^hT^*(V)$ the  $h^{th}$  symetric power of the bundle  $T^*(V)$ of 1-forms, we call   {\it  principal  symbol} of $D_h$
		the restriction $$\sigma_h(D):S^hT^*(V)\otimes E\to S^{h-k}T^*(V)\otimes F$$ of  $D_h$ to the kernel $S^hT^*(V)\otimes E$ of the projection $J^hE\to J^{h-1}E$, which takes values in the kernel  $S^{h-k}T^*(V)\otimes F$ of the projection $J^{h-k}F\to J^{h-k-1}F$.
		
	Let $p$ be the rank of $E$ and assume the rank   of the morphism $(J^kE)_m\to F_m$ to be constant equal to $q$ (not depending on the point $m\in V$) : thus, replacing if necessary the initial $F$ by $Im(D)$, we may assume the vector  bundle  $F$ to have rank $q$,   and $D$ to be an epimorphism.
		
		We denote by  $$c(n,h):=\frac{(n-1+h)!}{(n-1)!\ h!}$$the dimension of the  vector space of homogeneous polynomials of degree $h$ with $n$ variables, which is also the rank of $S^hT^*(V)$, or the number of multi-indices $I=(i_1,...,i_n)$  of partial derivatives of order $|I|=h$ with respect to $n$ variables, (where
		$|I|=i_1+i_2+...+i_n$).
		
		\pagebreak
		Recall that the  homogenization of polynomials by adding one more variable implies the formula $$\sum_{\ell=0}^h c(n,\ell) =c(n+1,h).$$
	 Hence, $J^hE$ has rank $p.c(n+1,h) $, and $J^{h-k}F$ has rank $q.c(n+1,h-k) $. 	
		
		\n {\bf Definition 1 :}
		{\it  The differential operator  $D$  is said to be  {\rm ordinary} if, for any $h$ $(h\geq k)$, the principal symbol $\sigma_h(D)$ has maximal rank $inf \bigl(q.c(n,h-k)\ ,\ p.c(n,h)\bigr) $.}

		Because of the proposition 1 below,  we shall see that, for $D$ to be ordinary,  we need    in fact  to check only  a  {\it  finite} number of conditions, at least when $p<q\leq p.c(n,k)$. In this case,  define   the  integer $h_0$ as being the  biggest integer $h$ ($h\geq k$) such that $p.c(n,h)\geq  q.c(n,h-k)$. This makes sense because of the 
		
			\n {\bf Lemma  1 : }{\it Let  $$\varphi(h):= \frac{c(n,h-k)}{c(n,h)}\ \ \ \Bigl(=\prod_{r=h-k+1}^h\frac{r}{n-1+r}\Bigr).$$ The function $\varphi$
			is strictly increasing from $\frac{1}{c(n,k)}$ to 1 when $h$ goes from $k$ to $+\infty$.}  
		
		\n {\it Proof :} In fact, $\varphi(h)-\varphi(h-1) =\frac{k(n-1)}{(n-1+h)(n-1+h-k)}.\prod_{r=h-k+1}^{h-1}\frac{r}{n-1+r} $. This  is a strictly positive number for $h\geq k+1$. Moreover, since  any function $r\mapsto \frac{r}{n-1+r}$ goes to 1 when $r$ goes to $+\infty$, the product $\prod_{r=h-k+1}^h\frac{r}{n-1+r}$  goes also to 1 (the number  $k$ of factors in this product being fixed).
		
	\medskip 
		
			\centerline{\includegraphics[scale=.55]{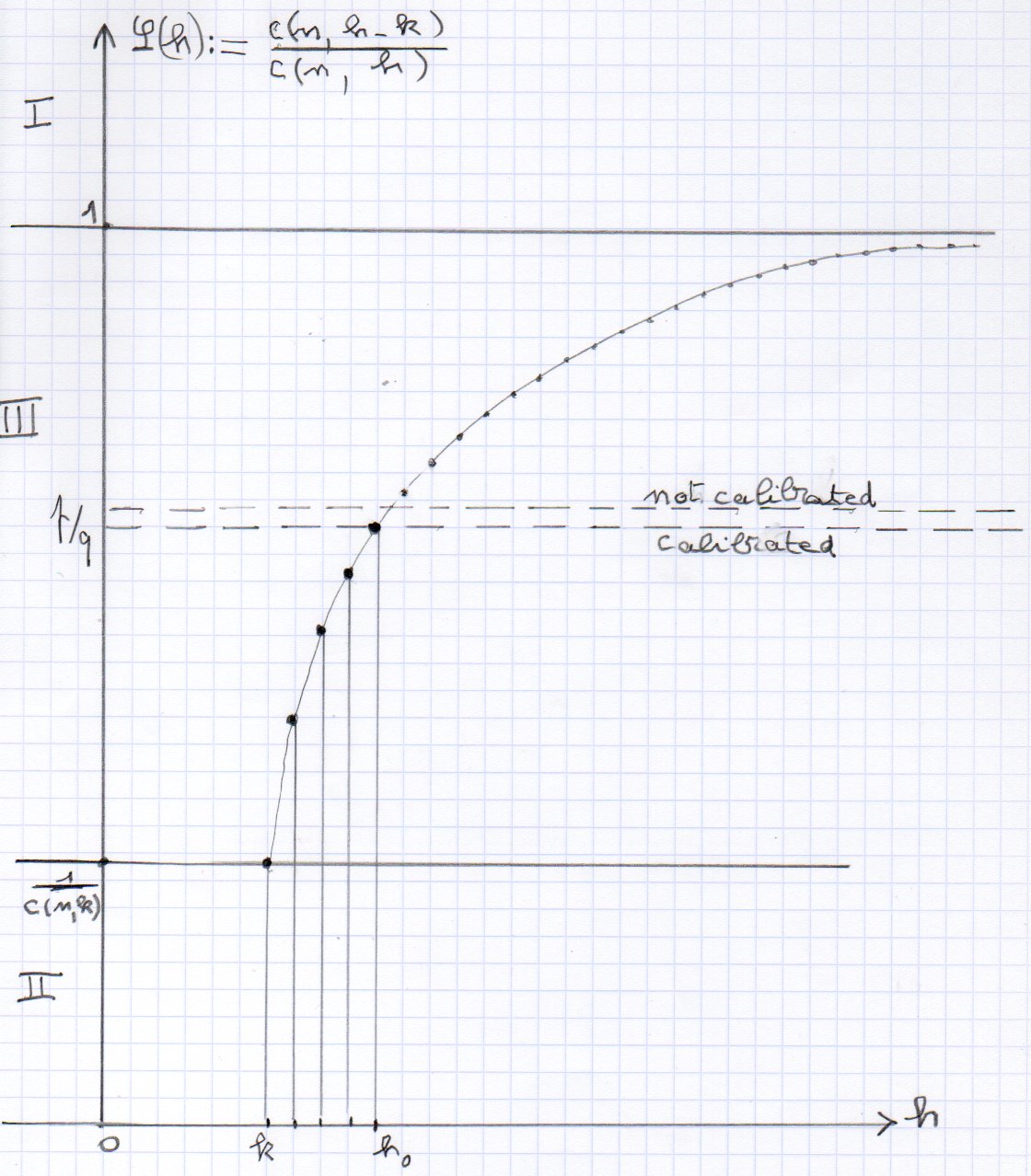}} 
		
	 \n {\bf Proposition 1 :} {\it 
	 	When $p<q\leq {p.c(n,k)},$  it is sufficient,  for $D$ to be ordinary, 
	 	that  the principal symbols $\sigma_h(D)$ have their rank maximal  for $k\leq h\leq h_0+1$ in the general case, $($and even for $k\leq h\leq h_0$ if,  moreover, $ p.c(n,h_0)=q.c(n,h_0-k). $$)$.} 
 	
 	\n The proof will be given in section 5.

		\n {\bf Definition 2 :}
		{\it  The differential operator  $D$  is said to be  {\rm calibrated} if it is ordinary, and  
			if  $$ p.c(n,h_0)=q.c(n,h_0-k). $$}

		\n {\bf Example : }Assume $n=2$, $p=2$,  $q=3$, and $k=1$. Then, $h_0=2$.  Locally,  given    local coordinates $(x,y)$ on the base space,  and after local trivializations of $E$ and $F$, a section of $E$ is defined by a function $s:(x,y)\mapsto s(x,y)$
		which takes values in the space of the 2-column vectors. The differential operator   writes $$<D, j^1s>=<A,s>+<B,s'_x>+ <C,s'_y>, $$
		where $A,B$ and $C$ denote matrices of size $3\times 2$ whose coefficients are scalar functions of $(x,y)$. Its first prolongation $D_2$ writes 
		$$<D_2, j^2s>=$$ $$<\begin{pmatrix}A\\A'_x\\A'_y\end{pmatrix},s>+<\begin{pmatrix}B\\A+B'_x\\B'_y\end{pmatrix},s'_x>+ <\begin{pmatrix}C\\C'_x\\A+C'_y\end{pmatrix},s'_y>+<\begin{pmatrix}0\\B\\0 \end{pmatrix},s''_{x^2}>+<\begin{pmatrix}0\\C\\B\end{pmatrix},s''_{xy}>+<\begin{pmatrix}0\\0\\C\end{pmatrix},s''_{y^2}>. $$
		Principal symbols write :
		$$\sigma_1(D)=\begin{pmatrix} B&C\end{pmatrix}\hbox{ and  } \sigma_2(D)=\begin{pmatrix} B&C&0\\0&B&C\end{pmatrix}.  $$
		Then $D$ is ordinary iff $\sigma_1(D)$ has rank $3$ and $\sigma_2(D)$ has rank $6$ at any point $(x,y)$, and is  then   calibrated  (since $\sigma_2(D)$ is   a square matrix $6\times 6$).
		
 We shall see other examples  in sections 6, when looking for abelian relations of  
		curvilinear webs.  
	
			\section{Claim of the main results}

		Let $D: J^kE\to F $ be a homogeneous linear differential operator of order $k$, with
		$n$, $p$, and $q$ as in section 1. 
				Associated to it is the  partial order differential  equation (shortly PDE) )   $$(*) \hskip 1cm <D,j^ks>  \equiv 0  . $$ 
				The kernel $R_h$ of $D_h:J^hE\to J^{h-k} F$ ($h\geq k$) is the set of the formal solutions of $(*)$ at order $h$.


 If  $D$ is  ordinary and if  $p<q\leq {p.c(n,k)},$ we shall prove ({\bf  theorem 1}) that      the dimension of the vector space   ${\mathcal S}_m$ of germs of solutions of $(*)$ at a generic point $m$ of $V$ is upper-bounded by the number  
$$ \pi(n,k,p,q):=p.c(n+1,h_0)-q.c(n+1,h_0-k) .$$
If $D$ is not calibrated, we need   to assume the framework to be analytic (real or complex).

 If    $D$ is moreover  calibrated (and then, we don't need analyticity), ${\mathcal E}:=R_{h_0-1}$ is a vector bundle\footnote{with the convention $R_{k-1}=J^{k-1}E$.}, whose rank is exactly $ \pi(n,k,p,q)$, and we define a tautological connection on it, having the following property : the sections $s$ of $E$ which are solutions of $(*)$ correspond bijectively to the sections $\s$ of $\mathcal E$ whose covariant derivative $\nabla \s$ vanishes by the map $s\mapsto 
j^{h_0-1}s$.    Hence, this  proves ({\bf  theorem 2}) that 
the dimension of   ${\mathcal S}_m$ is maximal iff the curvature of this connection vanishes. We prove also ({\bf  theorem 3}, called theorem of ``concentration") 
that, relatively to some convenient trivialization of ${\mathcal E}$, the only coefficients of the  curvature matrix  which may not vanish belong to the $p.c(n,h_0-1)-q.c(n,h_0-k-1)$ last lines.  

We then apply these results to the rank of curvilinear webs (results announced in [L]) :  for any  $d$-web by curves on a $n$-dimensional manifold, we can define  a differential operator of order $k=1$,
with $p=d-n$, $q=d-1$, which is always ordinary and calibrated, and  for which the solutions of $(*)$ are the $(n-1)$-abelian relations.  Using theorem 1, we recover   the Damiano's upper-bound  $$\sum_{h=0}^{d-n-1}\begin{pmatrix}n-2+h\\h\end{pmatrix}(d-n-h) \hskip 1cm	\Biggl(= \begin{pmatrix}d-1\\ n\end{pmatrix}\Biggr)$$
for the rank of the web $([D_1],[D_2]$). Using  theorem 2, we see that the Damiano's bound is reached iff  the curvature of the web vanishes ({\bf  theorem 4}). Moreover, after theorem 3,  the whole curvature -when looked as a matrix whose coefficients are 2-forms- is concentrated in the $\begin{pmatrix}d-3\\ n-2\end{pmatrix}$
last lines of the matrix  for a convenient trivialization. 

\section{ An upper bound for the dimension of the space of germs of solutions}

We want to study the tower 
$$R_\infty\to ...\to R_h\buildrel{\pi_h}\over\longrightarrow R_{h-1}\to ...\to R_{k=1}\to R_k\to J^{k-1}E .$$
Begining by the step $R_h\to R_{h-1}$ ($h\geq k$), let $r_{0}\in R_{h-1}$ and $\Sigma(r_{0})$ be 
the linear system of equations	$$ <D_h,r>=0,\ \  \ \pi_h(r)=r_{0}$$
whose solutions are the elements $r$ of $R_h$ projecting onto $r_0$. The homogeneous associated linear  system writes $<\sigma_h(D), r>=0$, where the principal symbol $\sigma_h(D):S^hT^*(V)\otimes E\to S^{h-k}T^*(V)\otimes F$ is   of maximal rank $inf\bigl(p.c(n,h)\ ,\ q.c(n,h-k)\bigr)$ if $D $ is assumed to be ordinary. Hence,  
$\Sigma(r_{0})$ is a system of $q.c(n,h-k)$ linearly independant scalar equations, with $p.c(n,h)$ scalar unknowns. The set of solutions is then an affine space of dimension $p.c(n,h)- q.c(n,h-k)$ or is over-determined according that $\frac{p}{q}$ is bigger or strictly smaller than $\varphi(h)= \frac{c(n,h-k)}{c(n,h)}$. After the snake's lemma, the projection $R^k\to J^{k-1}E$ is surjective when $D$ is ordinary and when  $q\leq p.c(n,k)$. In fact, under these assumptions, the cokernel of the symbol $\s_k(D)$ in the diagram below vanishes, while  the sequence $R_k\to J^{k-1}E\to Coker\bigl(\s_k(D)\bigr)$ is exact. 
$$\begin{matrix}
	
	&0&\to&Ker\bigl(\s_k(D)\bigr)&\to&R_k&\to&J^{k-1}E&---& ---\\
	
	&&&\downarrow&&\downarrow &&\hskip 2mm\downarrow =&&\hskip 1.15cm |\\
	
	&0&\to&S^kT^*V\otimes E&\to&J^kE&\to&J^{k-1}E&\to&0\hskip 1cm |\\

	&---&---&--\hskip .4cm\downarrow\s_k(D)&---&\hskip .3cm\downarrow D&---&\downarrow&-\!-\!-&--- \\
	
	
	|&0      &\to&F&\buildrel{\cong}\over\to&F&\to &0&\to&\hskip -1cm 0 \\
	
	|&&&\downarrow&&\downarrow\\
	
	&---&--&\hskip -1.1cm --\to\ Coker\bigl(\s_k(D)\bigr)&\to &0\\ 
	
\end{matrix}$$
Since $\varphi(h)$ increases from $\frac{1}{c(n,k)}$ to 1 when $h$ goes from $k$ to $+\infty$, we distinguish three ranges:

\n (I)- If $q\leq p$, there are very few independant equations with respect to the number of scalar variables, the inequality  $\frac{p}{q}> \varphi(h)$ holds  whatever be $h$ ($h\geq k$), and  all $R_h$ are vector bundles of rank $$\rho_h=p.c(n+1,k-1)+\sum_{\ell=k}^h\bigl(p.c(n,\ell)- q.c(n,\ell-k)\bigr) \ \ \Bigl(=p.c(n+1,h)-q.c(n+1,h-k)\Bigr).$$ 
The dimension $\rho_\infty$ of formal solutions     at a point $m$ of $V$ is $+\infty$, and $dim ({\mathcal S}_m)$ may be anything, including $+\infty$.

\n (II)- If $q>p.c(n,k)$, there are   many  independant equations with respect to the number of variables,
$\frac{p}{q}< \varphi(h)$   whatever be $h$ ($h\geq k$),  so that  $dim({R_\infty})_m$ (and a fortiori $dim ({\mathcal S}_m)$ if the framework is analytic) may not be bigger than the rank
$p.c(n+1,k-1)$ of $J^{k-1} E$.

\n (III)-  If $p<q\leq p.c(n,k)$, all $R_h$ such that $k\leq h\leq h_0$ are vector bundles of rank $$\rho_h=p.c(n+1,k-1)+\sum_{\ell=k}^h\bigl(p.c(n,\ell)- q.c(n,\ell-k)\bigr) \ \ \Bigl(=p.c(n+1,h)-q.c(n+1,h-k)\Bigr).$$ For $h>h_0$, all systems $\Sigma(r_{0})$ are over-determined, and have at most 1 solution ; thus, the dimension of formal solutions  at a point $m$ of $V$ is at most equal to $\rho_{h_0}$. If the framework is analytic (real or complex), $dim ({\mathcal S}_m)$ is a fortiori upper bounded by this number. Hence

\n {\bf Theorem 1 :}
{\it Assume  the differential operator $D$ to be  ordinary,  the framework to be analytic\footnote{After theorem 2 below, we shall see that we don't need in fact the assumption of analyticity if $D$ is calibrated.} (real or complex)
	and  $p< q\leq p.c(n,k) $ $(h_0\geq k)$. Then, the space ${\mathcal S}_m$ of germs of  solutions of the partial order equation $(*)$ at any point $m$ of $V$ is  a finite dimensional vector space of dimension at most equal to $\pi(n,k,p,q):=\rho_{h_0}$
	$$\pi(n,k,p,q):= p.c(n+1,k-1))+\sum_{h=k}^{h_0}\begin{pmatrix}n-1+h\\h\end{pmatrix}.\bigl(p-q.\varphi(h)\bigr) ,\hskip 1cm \biggl(=p.c(n+1,h_0)-q.c(n+1,h_0-k)\biggr).$$
	Moreover the germ of a solution at the point $m$ is completely defined by its $h_0$-jet at this point.}

\n Notice that  the summation may be done only up to $h_0-1$ in the calibrated case, since then  the last term vanishes.

\section{Calibrated differential operators, and connections}

We assume  the differential operator $D$ to be ordinary, and such that  $p< q\leq p.c(n,k) $. 
If $D$  is moreover calibrated, then $\frac{p}{q}=\frac{c(n,h_0-k)}{c(n,h_0)} $ : all   systems $\Sigma_{r_0}$ are cramerian  for $r_0\in R_{h_0-1}$,  and the natural  projection $R_{h_0}\buildrel{\pi_0}\over \longrightarrow R_{h_0-1}$ is now an isomorphism. Recalling that 
$R_{h_0}$ is the intersection of $J^1R_{h_0-1}$ and $J^{h_0}E$ in $J^1J^{h_0-1}E$, let $\iota $ be the inclusion of $R_{h_0}$ into $J^1(R_{h_0-1})$ ; the composed morphism $u:=\iota\scirc (\pi_0)^{-1}$ is a splitting of the projection $J^1(R_{h_0-1})\to R_{h_0-1}$ and defines therefore a connection on the vector bundle ${\mathcal E}:=R_{h_0-1}$ : for any section $\sigma$ of $\mathcal E$,  the covariant derivative is given by the formula$$\nabla \sigma =j^1\sigma-u(\sigma),$$ this difference being in the kernel $T^*(V)\otimes {\mathcal E}$  of the projection $J^1{\mathcal E}\to {\mathcal E}$ :
$$0\to T^*(V)\otimes {\mathcal E}\buildrel{\buildrel{\nabla}\over \longleftarrow}\over\longrightarrow J^1{\mathcal E}\buildrel{\buildrel{u}\over \longleftarrow}\over\longrightarrow  {\mathcal E}\to 0$$ 

\n {\bf Theorem 2 :} {\it If $p< q\leq p.c(n,k) $, and if  $D$ is  ordinary and  
	calibrated, then :
	
	$(i)$ the solutions of $(*)$ are the sections  $s$ of ${\mathcal E}$ such that $\nabla (j^{h_0-1}s)=0$,
	
	$(ii)$  the dimension of ${\mathcal S}_m$ is at most equal to $\pi(k,n,p,q)$, and is equal to $\pi(k,n,p,q)$  iff the curvature of the previous connection vanishes.
} 

\n {\it Proof  :}

Let $s$ be a section of $E$. It is a solution of $(*)$ iff $j^{h_0}s$ belongs to $R_{h_0}$. Since $R_{h_0}$ is the intersection $J^1R_{h_0-1}\cap J^{h_0}E$ in $J^1J^{h_0-1}E$, and since we know already that $j^{h_0}s$ is always a section of $J^{h_0}E$, it is equivalent to say that it is also a section of $J^1R_{h_0-1}$ which may be written  $j^1(j^{h_0-1}s)=\iota((\pi_0)^{-1}(j^{h_0-1}s))$, i.e. $\nabla (j^{h_0-1}s)=0$. Therefore,  the equations $\nabla j^{h_0-1}s=0$  and $D_{h_0}j^{h_0}s=0$
have the same solutions : this proves part $(i)$ of the theorem.

The curvature vanishes, iff the module of sections of $\mathcal E$ has a  local basis made of  $\nabla$-invariant sections. Since the rank of $\mathcal E$ is equal to the maximal possible dimension of ${\mathcal S}_m$, this proves  part $(ii)$. 

\rightline{QED}
\n {\bf Remarks :} Notice that the construction of the connection  doesn't require assumption of analyticity. 
Therefore, in the calibrated case, the upper-bound given in theorem 1 is still available in the real differentiable case without assumption of analiticity. 

Moreover,  in the calibrated case,  we  don't need to prove that $\s_h(D)$ has maximal rank for $h>h_0$, since   the solutions are the $\nabla$-invariant sections of a bundle $\mathcal E$,  which has the required rank $\pi(k,n,p,q)$. 

\n {\bf Concentration of the curvature :}

We first define a filtration\footnote{There is no reason a priori for this filtration to be preserved by the connexion, since $\pi_0$ has nothing to see with it.}  of   ${\mathcal E}\ \bigl(=R_{h_0-1}\bigr)$ 
$$F_{k-1}({\mathcal E})={\mathcal E}\supset F_{k}({\mathcal E})\supset ...\supset F_{h}({\mathcal E}) \supset F_{h+1}({\mathcal E})\supset ...\supset F_{h_0-1}({\mathcal E})\supset F_{h_0}({\mathcal E})=0$$
by setting 
$F_h({\mathcal E}):=Ker({\mathcal E}\to R_{h-1})$ (with $k-1\leq h\leq h_0$, and the convention $R_{k-1}=J^{k-1}E$).
This filtration is induced by the similar  filtration $F_h(J^{h_0-1}E):=Ker(J^{h_0-1}E\to J^{h-1}E)$ of $J^{h_0-1}E$.

 Let  $Gr({\mathcal E}):= \bigoplus_{h=k-1}^{h_0-1}\bigl(F_h({\mathcal E})/F_{h+1}({\mathcal E})\bigr)$ be the 
 graded vector bundle associated to this filtration, 
  which is a sub graded vector bundle of  the graded vector bundle $Gr(J^{h_0-1}E)$  associated to the filtration of $J^{h_0-1}E$. 
   
We may represent a section $s$ of $E$ as a $p$-column vector (whose componants are scalar functions) once $E$ has been  locally  trivialized, and a local section $s_h$ of $Gr^h(J^{h_0-1}E)$ as a family $(s_I)_{|I|=h}$ indexed by by the set of multi-indices $I$  of order $h$, once given a system $(x_i) _i$ of local coordinates on $V$ with respect to which $I$-derivations may be done. 

Let $\rho:{\mathcal E}\to Gr({\mathcal E}$ be the natural projection (isomorphism of non-graded vector bundle).

A  local trivialization of ${\mathcal E}$ is said to be ``adapted'') if  there exists a local trivialization of $E$ and a system of local coordinates $(x_i)_i$ such that the  following condition holds : 

\n {\it  for each $h$, ($k\leq h\leq h_0-1$), there exists a subset ${\mathcal I}_h$ of $\frac{1}{p}\bigl(p.c(n,h)-q.c(n,h-k)\bigr)$ multi-indices $I$ among those of order $|I|=h$ $($with the convention that ${\mathcal I}_{k-1} $ is the set of all $c(n+1,k-1)$ multi-indices of order $\leq k-1)$, such that the family $\{s_I\}_{I\in {\mathcal I}_h, \  k-1\leq h\leq h_0-1}$ be the trivialization of 
$Gr^h({\mathcal E})$ corresponding by $\rho$ to the given trivialization of  ${\mathcal E}$}.

 
\n {\bf Theorem 3 :} {\it Let us represent   the curvature   of the above tautological connection $\nabla$, relatively to some trivialization of ${\mathcal E}$,   as  a 2-form $K$ on $V$ with coefficients in the space of the  $\pi(k,n,p,q)\times \pi(k,n,p,q)$ matrices.\hb
	\indent 	If this trivialization is adapted, the only lines of these matricial coefficients which may have   non-vanishing terms  are among  the $\bigl(p.c(n,h_0-1)-q.c(n,h_0-1-k)\bigr)$ last ones.  Equivalently, the lines from 1 to  \break    $p.c(n+1,h_0-2)-q.c(n+1,h_0-2-k)$  are lines of zeros.}

\n The proof will be given in section 5.

\n {\bf Particular case  $k=1,\ h_0=1$ :}

Then, if $D$ is ordinary,  the space 
${\mathcal S}_m$ of germs of solutions of $(*)$ has maximal dimension $p$. If $D$ is  moreover calibrated (i.e. $q=pn$), $\sigma_1(D)$ is invertible : if the $p$-column-vector  $f:=(f_u)_{u=1,...,p}$  denotes the componants of a section of $E$ with respect to some local trivialization, the equation $(*)$ means  that the partial derivatives $ f'_i:=\frac{\partial f}{\partial x_i}$  of a solution with respect to some local coordinates $(x_i)_i$ must be some linear  function of $f$ : $ \partial_i f =<A_i,f>$ for some $p\times p$-matrix $A_i$, and the vanishing of the curvature of the tautological connection defined above means the commutativity of the second partial derivatives of the  : $$\Bigl(<A_j,f>\Bigr)'_i-\Bigl(<A_i,f>\Bigr)'_j=0 \hbox{ for any } i<j .$$

For instance, if $p=1$,  $q=n$, $f$ and   $A_i$ are scalar functions and the non-vanishing solutions of $(*) $ are defined by the solutions of the equation $\frac{df}{f} =\sum_i A_i\ dx_i$ ; the corresponding connection writes locally $$\nabla_if=df-\omega.f \hbox{ \hskip 1cm where }\omega:=\sum_i A_i\ dx_i.$$There exists or not such solutions    according to the fact the 1-form $\omega$ is closed or not,  the curvature of the conection being equal to $ d\omega$ up to sign.

\section{Proofs :}

\subsection{Generalities}
\n Let $(\partial_i)_{1\leq i\leq n}$ be a local basis of the module of vector fields\footnote{It will be convenient for the sequel to accept that these vector fields may not commute.} on an open set  $U$ of the ambiant manifiold $V$. 
Let $(f_v)_{1\leq v\leq p}$   be the componants of a section of $E$  with respect to some local trivialization $(\epsilon _v)_{1\leq v\leq p}$ of $E$, and 
$(\tau _u)_{1\leq u\leq q}$ be some some local trivialization of $F$. 

\subsubsection{Indices of derivations :}
First, if $I=(i_1,i_2,...,i_n)$ denotes a multi-index of height  $|I|=h$, (where $|I|:=i_1+i_2+...+i_n$), the derivative of higher order $g'_I$ of a funcion $g$ with scalar or vectorial  values will mean that we have first taken $i_1$ times the derivative with respect to $\partial_1$, then $i_2$ times  with respect to $\partial_2$, etc....Since  $\partial_i.(\partial_j.g)=\partial_j.(\partial_i.g)+[\partial_i,\partial_j].g$, any sequence of  derivatives made in a different order may be replaced by the previous one, modulo derivatives of strictly smaller order. In particular, this change of    ordering would  not affect the description of the symbols $\s_h(D)$  which are the terms of highest order in the expression of $D_h$.

We
need an  ordering for  all $n$-multi-indices of height $\leq h_0$ using  a bijection $$LL:\{0,1,\cdots ,c(n+1,h_0)-1\}\to \{multi-indices\},$$ that we choose so that 

$LL(0)=(0,0,...,0)$, 

$LL(i)$ is equal, for any ($1\leq i\leq n$), to  the   multi-index $1_i:=(0,0,...,0,1,0,...,0)$ with 1 at the  $i^{th}$ place,

and, more generally, any multi-index  of heigt $h$ receives  a place between $c(n+1,h-1)$ and $c(n+1,h)-1$.

\n The rule used in the  Maple  programm of [L] : $(i_1,i_2,...,i_n)<(j_1,j_2,...,j_n)$ iff 

- either   $i_1+i_2+...+i_n<j_1+j_2+...+j_n$,

- or $i_1+i_2+...+i_n=j_1+j_2+...+j_n$, and   $i_k<j_k$, where $k$ denotes the highest index $h\in \{1,...,n\}$ such that $i_h\neq j_h$ : $ (i_1,i_2,...,i_k,a_{k+1},...,a_n)<(j_1,j_2,...,j_k,a_{k+1},...,a_n) $.

\n For example, if $n=3$, we get the order :\hb 
$(000), (100, 010,001)  , (200,110,020,101,011,002), (300,210,120,030,201,111,021,102,012,003), $\hb 
$(400,310,220,130,040,301,211,121,031,202,112,022,103,023,004), ...$.

\n If $H=LL(t)$, we shall write sometimes $g'_t$ instead of $g'_H$. Observe that, for $t=1,...,n$, there is no ambiguity since then $LL(t)=1_t$. We shall also write $|t|$ instead of $|LL(t)|$ (if $1\leq t\leq n$ for instance, $|t|=1$).
\subsubsection{Matrices $P_h$, $P_h$ and $Q_h$ :}
Then $D$ is locally defined by $q$ equations $( EQ_u)$, $u=1,...,q$ :
$$( EQ_u)\hskip 1cm \sum_{v=1}^p\Bigl(\sum_{K,|K|=k}  <M_{u}^{K,v},(f_v)'_K>+\sum_{K,|K|<k} <  M_{u}^{K,v},(f_v)'_K>  \Bigr)\equiv 0,$$
where the coefficients $M_{u}^{K,v}$,  are known, and the scalar functions $f_v$ unknown. 

\n We can write shortly these $q$ equations together as $<M_1,j^1f>=0$, or more precisely
$$( \Sigma_k)\hskip 1cm \sum_{k,|K|=k}^n  <T_K,f'_K>+< Q_k,j^{k-1}f >\equiv 0,$$
where  $f$ denotes the  $p$-column vectors  of  componants $(f_u)_u$, $T_K$ the matrix $(\!(M_{u}^{K,v})\!)_{u,v}$ of size $q\times p$ with $|K|=k$, and $Q_k$ some known matrix of size $q\times p.c(n+1,k-1)$ defined by the coefficients 
$M_{u}^{K,v}$ for $|K|<k$.

\n The   principal symbol $\sigma_k(D)$ is then  represented  by the  matrix  of size  $q\times p.c(n,k)$.   $$P_k=(T_{K_1}\ T_{K_2}\    ...\ T_{K_{c(n,k)}}), $$ and $R_k$ is defined by the kernel of the $q\times c(n+1,k)$-matrix $M_k=(Q_k\ P_k)$.

\n The equation $D_h(j^hs)=0$ is then equivalent to the set of the $q.c(n+1,h-k) $ equations $( EQ_u)'_H$  ($|H|\leq h-k$) and  $ p. c(n+1,h)$ unknowns  $((f_v)'_K)_{|K|\leq h}$, written 
shortly $<M_h,j^h f>\equiv 0$ with $$M_h=\begin{pmatrix}M_{h-1}&0\\
	Q_{h}&P_{h}\
\end{pmatrix} .$$
being such that equation $( EQ_u)'_H$ writes :
$$\sum_{K,v}M_{H,u}^{K,v}.(f_{v})'_K\equiv 0 \hbox{ with } |K|\leq |H|+k \hbox{ and  } 1\leq v\leq p.$$
The  coefficients $M_{H,u}^{K,v}$ of  $M_h$ (hence those of    $P_h$ and  $Q_h$) are gotten inductively.   

\n Then, if $H =LL(t)$  and $K=LL(s)$ ($|K|\leq h$), the coefficient $M_{H,u}^{K,v}$ of $M_h$  will be on the  line $i=u+qt$
and in the colum  $j=v+ps$. We write :
$$ M_{H,u}^{K,v}=F(u+q.t,v+p.s).  $$

\n Conversely, we can  recover  $(u,H)$ from $i$ and  $(v,K)$ from $j$, with the  maps 

\n $i\to H(i) =LL([\frac{i+q-1}{q}])$  and   $i\to u(i) =i-[\frac{i+q-1}{q}]$ on one hand,  \hb    $j\to K(j) =LL([\frac{j+p-1}{p}])$  and   $j\to v(j) =j-[\frac{j+p-1}{p}]$ on the other hand   : 
$$F(i,j) =M_{H(i),u(i)}^{K(j),v(j)}$$ 

\n For any index $LL(t)$ and for any integer $i=1,...,n$, we shall denote by $ad_it$ the integer such that the index of  derivation by  $LL(ad_it)$ will mean that we shall take the derivative one more time by $\partial_i$ than with $LL(t)$ : if we write  $LL(t)=(L(t)_1,...,L(t)_n)$ the componants of $LL(t)$, then  $L(ad_it)_j=L(t)_j$ if $j\neq i$, and $L(ad_it)_i=L(t)_i+1$. Hence, 
$(EQ_u)'_{ad_i(t)}$  is the derivative of the identity  $(EQ_u)'_t$ with respect to 
$\partial_i$.

\n We know already the initial equations corresponding to $t=0$ :
$$( EQ_u)\hskip 1cm \sum_{v=1}^p\sum_{s,|s|\leq k}<  F(u,v+p.s),(f_v)'_s> \equiv 0,$$
The other coefficients are computed  inductively according to the height  $|t|$ of $LL(t)$ : assume already known the coefficients 
$F(u+q.t, v+p.s) $ of the identity
$$(EQ_u)'_t\hskip 1cm \sum_v\sum_{s, |s|\leq |t|+k}F\bigl(u+q.t,v+p.s\bigr).(f_v)'_s\equiv 0.$$Then, by derivation with respect to $\partial_i$, 
$$ 
(EQ_u)'_{ad_it}\hskip 1cm\sum_v \sum_{s,|s|\leq |t|+k}\biggl(\Bigl(F\bigl(u+q.t,v+p.s\bigr)\Bigr)'_i.(f_v)'_s+
F\bigl(u+q.t,v+p.s\bigr).(f_v)'_{ad_is}\biggr)\equiv 0,
$$hence 
$$F(u+q.ad_i(t),v+p.s) =0 \hbox{\hskip 3cm for $|s| >|t|+k+1$} ,$$
and the following formulae hold :

\n {\bf Proposition 2 :}

$F\bigl(u+q.ad_it,v\bigr)=\Bigl(F\bigl(u+qt,v\bigr)\Bigr)'_i$ \hskip 4.7cm for $s=0$,

$F\bigl(u+q.ad_it,v+p.s\bigr)=\Bigl(F\bigl(u+qt,v+p.s\bigr)\Bigr)'_i$ \hskip 3cm for  $1\leq |s|\leq h(t)+k$ and  $s\notin Im(ad_i)$,

$F\bigl(u+q.ad_it,v+p.ad_ir\bigr)=\bigl(F(u+qt,v+p.ad_ir)\bigr)'_i+F\bigl(u+qt,v+p.r\bigr)$   \hskip 1.3cm for $s=ad_ir$ and $ |s|\leq h(t)+k$,

$F\bigl(u+q.ad_it,v+p.ad_ir\bigr)=F\bigl(u+qt,v+p.r\bigr)$\hskip 3.2cm  for $s=ad_ir$ and  $ |s|= |t|+k+1$,

$F\bigl(u+q.ad_it,v+p.s\bigr)=0$ \hskip 6.1cm for $ |s|= |t|+k+1$ and   $s\notin Im(ad_i)$.

 \subsection{Proof of proposition 1}
 
  Denote by $K_t^s$ the  $q\times p$ bloc $(\!(F(u+q.t,v+p.s))\!)_{u,v}$.
In particular, 
 $P_h$ ($h\geq k$)  may be decomposed  into  $c(n,h-k).c(n,h)$ blocs $K_t^s$ of  size  $q\times p$, corresponding to the case $|t|=h-k$ and $|s|=h$, and defined inductively from $P_k=(T_{K_1}\ T_{K_2}\    ...\ T_{K_{c(n,k)}}), $  by\hb 
 \indent $K_t^{ad_ir}=K_t^r$ if $s=ad_ir$, $|r|=h-1$, $|t|=h-k$, \hb  
 \indent  and $K_t^s=0$ if $s\notin Im(ad_i)$, $|s|=h$, $|t|=h-k$. 
 
 \n Hence, by induction, we observe that the set of all blocs $\{K_t^s\}_s$ for a fixed $t$, with   $ |s|= |t|+k$  contains  all blocs  $T_{K_1}\ T_{K_2}\    ...\ T_{K_{c(n,k)}}$, the remaining blocs being   blocs of zeros.

 \n Example for  $n=3$, $k=2$   
 (it  is understood that blank  blocs are blocs of  zeros) :
 
 $P_2=\left(\begin{array}{cccccc}
 	T_{K_1}& T_{K_2}	& T_{K_3}&	T_{K_4}& T_{K_5}	& T_{K_6}\\
 \end{array}\right)$
 
 $P_3=	\left(\begin{array}{cccccccccc}
 	
 	T_{K_1}&T_{K_2}&T_{K_3}&&	T_{K_4}& T_{K_5}&	& T_{K_6}&&\\
 	
 	&T_{K_1}&T_{K_2}&T_{K_3}&&	T_{K_4}& T_{K_5}&	& T_{K_6}&\\
 	
 	&&&&T_{K_1}&T_{K_2}&T_{K_3}&	T_{K_4}& T_{K_5}&	 T_{K_6}\\	
 \end{array}\right)$
 
 
 		
 		
 		
 		
 		
 		

 		\n The  proof goes by induction on $h$. 
 		Assume that all column-vectors of $P_h$ are linearly independant for some $h\geq h_0$. Then we shall prove the same for all column-vectors of $P_{h+1}$. For that, it is sufficient to prove that their projection onto some sub-space gotten by suppressing  some componants are linearly independants :
 		since the set of all blocs $\{K_t^s\}_s$ for a fixed $t$, with   $ |s|= |t|+k$  contains  all blocs  $T_{K_1}\ T_{K_2}\    ...\ T_{K_{c(n,k)}}$, all  column-vectors of $P_{h+1}$ corresponding to blocs $K_t^s$ such that $s=ad_1t$ are linearly independant ; this is still true if we add the column-vectors of $P_{h+1}$ such that  $s=ad_2t$ which did not appear erlier ; this is still true if we add the column-vectors of $P_{h+1}$ such that  $s=ad_3t$ which did not appear erlier, etc... . Finally all column-vectors of $P_{h+1}$ are linearly independant. Hence, the rank of $P_h$ is maximal.  According to the assumption of the proposition 1, and according to the definition of $h_0$,   all column-vectors  of $P_{h_0+1}$ in the general case (of $P_{h_0} $ in the calibrated case) are linearly independant.
 		
 		\rightline{QED}
 
  \subsection{Proof of theorem 3}

  Let $(s_I )_I$ be the components of a section $s$ of ${\mathcal E}$ relatively to any adapted trivialization of ${\mathcal E}$, using   $(x_1,...,x_n)$ as local coordinates on $V$ ;   writing shortly $\partial_i$ unstead of $\frac{\partial }{\partial x_i}$  and $\nabla _i$ unstead of $\nabla _{\partial_i} $,  we get : 
  $$\begin{matrix} \bigl(\nabla_i s\bigr)_I&=&\partial_i s_I&-&s_{I+1_i}&\hbox{ if } |I|\leq h_0-2,\\
  	&&&&&\\
  	&=&\partial_i s_I&-&<U,s>_{I+1_i}&\hbox{ if } |I|= h_0-1,\\\end{matrix}$$
  	where $U:=-(P_{h_0})^{-1}.Q_{h_0}$  from $R_{h_0-1}$ into $Ker(J^{h_0}E\to J^{h_0-1}E)$,\hb 
  	and $I+1_i$ denotes the multi-index of derivation gotten from $I$ by one more derivation with respect to $\partial_i$.

 \n  Therefore 
  $$\begin{matrix} 
  	\bigl(\nabla_i\nabla_j s\bigr)_I&=&\partial_i(\partial_j s_I-s_{I+1_j})&-&(\partial_j s_{I+1_i}-s_{I+1_j+1_i})&\hbox{ if } |I|\leq  h_0-3\\
  	&&&&&\\
  	&=&\partial_i(\partial_j s_I-s_{I+1_j})&-&(\partial_j s_{I+1_i}-<U,s>_{I+1_j+1_i})&\hbox{ if } |I|=  h_0-2.
  \end{matrix}$$
  In both cases, the right-hand term  is symmetric in $(i,j)$, hence the vanishing of the coefficients 
  $\bigl(\nabla_i\nabla_j s-\nabla_j\nabla_is\bigr)_I $ in  the curvature, when  $|I|\leq h_0-2$. 
  
  \rightline{QED}

   \section{Application to curvilinear webs}

The framework being still    real $C^\infty$ or holomorphic, we are given  
a $d$-web $\mathcal W$ by curves  on  a  $n$-dimensional manifold $V$ ($d>n$),    defined  by $d$ foliations   ${\mathcal F}_\lambda$ of codimension $n-1$. Denote  by $T_\lambda$ (resp. $N_\lambda$)  the  tangent  (resp. normal) bundle to  ${\mathcal F}_\lambda$ :
$$o\to T_\lambda\to T(V) \to N_\lambda\to 0.$$
Let  $\partial_\lambda$ denote  a non-vanishing vector field generating locally  $T_\lambda$ on  a neighborhood $U$ of a point $m$ of $V$.\hb We shall assume that any  $n$ of these vector fields among the $d$'s are  linearly  independant at any point of $ U$.  

The dual sequence  $$0\to N^*_\lambda\to T^*({V}) \to T^*_\lambda\to 0$$
allows to identify   $N^*_\lambda$ to the bundle of  1-forms $\eta$ on $V$ which are ${\mathcal F}_\lambda$-semi-basic (i.e. such that the interior product  $\iota_{\partial_\lambda}\eta$ vanishes). Similarly,   $\bigwedge^{*} N^*_\lambda$ is identified to the forms on $V$  which are  ${\mathcal F}_\lambda$-semi-basic . The map $(\eta_\lambda)_\lambda\to \sum_\lambda \eta_\lambda$ from $\bigoplus_\lambda \bigwedge_\lambda^{n-1}N^*_\lambda$ to $\Lambda^{n-1}T^*({V})$ has maximum rank  $n$.
The kernel  $E$ of this map is therefore a vector budle of rank  $d-n$ : 
$$0\to  E\to \bigoplus_{\lambda=1}^d \bigwedge^{n-1}N^*_{\lambda}\to \bigwedge^{n-1}T^*({V}).$$
A section of this bundle (i.e. a family $(\eta_\lambda)_\lambda$ such that  $\sum_\lambda\eta_\lambda=0$) is called a $(n-1)$-abelian relation of $W$ (or shortly \emph{ an abelian relation}), if for any $\lambda=1,...,d$, each  $\eta_\lambda$ is  ${\mathcal F}_\lambda$-basic (i.e. not only ${\mathcal F}_\lambda$-semi-basic, but also ${\mathcal F}_\lambda$-invariant : $L_{\partial_\lambda} \eta_\lambda=0$,  where 
$L_{\partial_\lambda}=\iota_{\partial_\lambda}\scirc d+d\scirc \iota_{\partial_\lambda}$ denotes the Lie derivative). Hence abelian relations are sections $s=(\eta_\lambda)_\lambda$ of $E$ solutions of the equation ${\mathcal D}s=0$, 
where  ${\mathcal D}\bigl((\eta_\lambda)_\lambda\bigr):=\bigl(L_{\partial_\lambda}\eta_\lambda\bigr)_\lambda$ is 
associated to a linear homogeneous differential operator of  order $k=1$
$$ J^1E\buildrel{D}\over\longrightarrow \Bigl(\wedge^{n-1}T^*({\mathcal U})\Bigr)^{\oplus d}.$$  
The set of  $n-1$-abelian relations has a natural    structure of vector space,  whose  dimension (may be infinite a priori) is called  the $(n-1)$-rank (or shortly the {\it rank}) of the web.

We proved in [L] the 

\n {\bf Proôsition 3 :}{\it The differential operator $D:J^1E\to \Bigl(\wedge^{n-1}T^*({\mathcal U})\Bigr)^{\oplus d}$,  just defined above for any $d$-web of curves in a $n$-dimensional ambiant manifold, is always ordinary and calibrated, with $p=d-n$, and $q=d-1$. }

Thus $ \frac{p(n-1)}{q-p}$ is equal to the integer $h_0=d-n.$
The inequalities $p<q\leq np$ holds  since $d \geq n+1 $.
Therefore, after theorem 1,  we recover the upper-bound  
$$\sum_{h=0}^{d-n-1} \begin{pmatrix}n-1+h\\h\end{pmatrix}.(d-n-h) \hskip 1cm \Biggl(=\begin{pmatrix}d-1\\n\end{pmatrix}\Biggr)$$given by Damiano ([D1]) for the rank of the web. On the other hand, since  $D$ is ordinary and calibrated, there is a tautological connection on ${\mathcal E}:=R_{d-n-1}$   and,   as a corollary of the  theorems 2 and 3, we get  :

\n {\bf Theorem 4 :}

{\it The Damiano's    bound for the rank of a curvilinear bound is reached iff the curvature of the  tautological connection on $R_{d-n-1}$ defined above vanishes.

\n Moreover, this  curvature is concentrated in the $\begin{pmatrix}d-3\\n-2\end{pmatrix}$ last lines
of the matrix, when an adapted trivialization is used.} 

In [L], we wrote a Maple program for computing the curvature of the following $(n+3)$-web $W_c$ defined,  relatively to some system  $(x_i)$ of affine   coordinates  on a convenient open set of a $n$-dimensional projective space,   by the vector fields 

$$\partial_i=\frac{\partial }{\partial x_i} \hbox{  for any $i=1,\cdots,n$}\ ,$$
$\partial_{n+1}= \frac{1}{x_n+c}\sum_i (x_i+c)\ \partial_i$,\hskip .4cm 
$\partial_{n+2}=\frac{1}{x_n-1-c}\sum_i (x_i-1-c)\ \partial_i$,\hskip .4cm 
$\partial_{n+3}=\frac{1}{x_n(x_n-1)}\sum_i x_i(x_i-1)\ \partial_i,$

\n where $c$ is a scalar parameter. For $c=0$,  we get the exceptional web $W_{0,n+3}$ (the Bol's web for $n=2$). 
It is known, since Bol ([Bo]) for $n=2$,   Damiano  ([D1],[D2]) for $n$ even and Pirio ([Pi3]) for $n$ odd) that $W_{0,n+3}$ has a maximal rank. We recovered this result,   at least\footnote{Theorically, the program works whatever be $n$, but the time  of computation   becomes too long with a small computer.}
for $n=2$ or $3$, by proving  that the curvature vanishes.

\n {\bf Remark :} 
  The trivialization used in [L] was not adapted (but had some other interest that we explained there). We wrote another    program with an adapted trivialization (which has 3 elements in the sections of $F_0({\mathcal E})$, 
  4 in those of $F_1({\mathcal E})$, and 3  in those of $F_2({\mathcal E})$,) and observed that, for $c\neq 0$, the curvature was effectively concentrated on the three last lines of the $10\times 10$ curvature matrix, as required by the concentration's theorem above.

\n {\bf Bibliography}

\vskip .5mm

\n [BB] W. Blaschke et G. Bol, {\it Geometrie der Gewebe}, Die Grundlehren der Mathematik 49, Springer, 1938.
\vskip .5mm

\n [Bo] G. Bol, {\it \"Uber ein bemerkenswertes F\"unfgewebe in der Ebene}, Abh. Math. Hamburg Univ., 11, 1936, 387-393.
\vskip .5mm



\n [CL] V. Cavalier, D. Lehmann, {\it Ordinary holomorphic webs of codimension one. } arXiv 0703596v2 [mathsDS], 2007, et Ann. Sc. Norm. Super. Pisa, cl. Sci (5), vol XI (2012), 197-214.
\vskip .5mm

\n [D1] D.B. Damiano, {\it Abelian equations and characteristic classes}, Thesis, Brown University, (1980) ;  American J. Math. 105-6, 
1983,  1325-1345.
\vskip .5mm

\n [D2] D.Damiano : {\it Webs  and characteristic forms on Grassmann manifolds}, Am. J. of Maths.105, 1983,  1325-1345.
\vskip .5mm

\n [DL1] J. P. Dufour, D. Lehmann, {\it Calcul explicite de la courbure des tissus calibr\'es ordinaires ; } arXiv 1408.3909v1 [mathsDG], 18/08/2014.
\vskip .5mm

\n [DL2] J. P. Dufour, D. Lehmann, {\it Rank of ordinary webs in codimension one : an effective method ; }\hb  arXiv 1703.03725v1 [math.DG], 10/03/2017. Pure and Applied Mathematics Quarterly, vol. 16, n. 5, 1587-1607, 2020. 
\vskip .5mm


\n [DL3] J.P. Dufour, D. Lehmann, {\it Etude des $(n+1)$-tissus de courbes en dimension $n$}, arXiv 2211.05188v1 [mathsDG], 09/11/2022, et 	Comptes Rendus Maths.  Ac. Sc. Paris, vol. 361, 1491-1497, 2023.
\vskip .5mm



\n [H1] A. Hénaut, {\it Planar web geometry through abelian relations  and connections}, 
Annals of Math. 159 (2004)  425-445.
\vskip .5mm

\n [H2]  A. Hénaut, {\it Formes différentielles abéliennes, bornes de Castelnuovo et géométrie des tissus}, Commentarii Math.  Helvetici, 79 (1), 2004, 25-57.
\vskip .5mm


\n [L]  D. Lehmann, {\it Courbure des tissus en courbes}, arXiv:2401.15988, v1(29/01/2024).
\vskip .5mm

\n [P] R.S.Palais. {\it Seminar on the Atiyah-Singer Index theorem}, Annals of Mathematics studies Princeton University press nr. 57.
\vskip .5mm

\n [Pa] A. Pantazi. {\it Sur la d\'etermination du rang d'un tissu plan.} C.R. Acad. Sc. Roumanie 4 (1940), 108-111.
\vskip .5mm

\n [Pi1] L. Pirio, {\it Equations Fonctionnelles Ab\'eliennes et G\'eom\'etrie des tissus},
Th\`ese de doctorat de l'Universit\'e Paris VI, 2004.
\vskip .5mm

\n [Pi2]  L. Pirio, {\it Sur les tissus planaires de rang maximal et le problème de Chern}, note aux C.R. Ac Sc. , sér. I, 339 (2004), 131-136.
\vskip .5mm

\n  [Pi3] L. Pirio : {\it On the $(n+3)$-webs  by rational curves induced by the forgetful maps on the moduli spaces ${\mathcal
		M}_{0,n+3}$}, arXiv 2204.04772.v1, [Math AG], 10-04-2022.
\vskip .2cm

\noindent Daniel Lehmann, former   professor  at the University of  Montpellier II.

\end{document}